\documentclass[dvips,ejs,preprint]{imsart}

\RequirePackage[OT1]{fontenc}
\RequirePackage{hyperref}

\usepackage[numbers]{natbib}
\usepackage{amsmath,amssymb,amscd,amsfonts,amsthm}

\RequirePackage{hypernat}

\doi{10.1214/11-EJS641}
\pubyear{2011}
\volume{5}
\issue{0}
\firstpage{1354}
\lastpage{1373}

\startlocaldefs

\usepackage{tikz}
\usetikzlibrary{shapes.geometric}
\numberwithin{equation}{section}
\theoremstyle{plain}

\theoremstyle{plain}
\newtheorem{theorem}{Theorem}[section]
\newtheorem{lemma}[theorem]{Lemma}
\newtheorem{proposition}[theorem]{Proposition}
\newtheorem{corollary}[theorem]{Corollary}
\theoremstyle{definition}
\newtheorem{definition}[theorem]{Definition}
\newtheorem{remark}[theorem]{Remark}

\def\denseset{W}
\def\indset{D}
\def\indI{\mbox{\tiny\rm I}}
\def\indJ{\mbox{\tiny\rm J}}
\def\indK{\mbox{\tiny\rm K}}

\def\indE{\mbox{\tiny\rm \indset}}

\def\ind#1{\mbox{\tiny\rm #1}}
\def\borel{\mathcal{B}}
\def\abstspace{\Omega}
\def\abstmeasure{\mathbb{P}}
\def\abstfield{\mathcal{A}}

\def\ie{i.e.\ }
\def\eg{e.g.\ }
\def\ind#1{\mbox{\tiny\rm #1}}

\def\xspace{\mathcal{X}}
\def\xspaceI{\xspace_{\indI}}
\def\xspaceJ{\xspace_{\indJ}}
\def\xspaceE{\xspace_{\indE}}
\def\borel{\mathcal{B}}
\def\borelI{\borel_{\indI}}

\def\borelE{\borel_{\indE}}
\def\top{\mathcal{T}}
\def\topI{\top_{\indI}}
\def\topJ{\top_{\indJ}}
\def\topE{\top_{\indE}}
\def\fJI{f_{\indJ\indI}}
\def\fKJ{f_{\indK\indJ}}
\def\fKI{f_{\indK\indI}}
\def\xI{{x_{\indI}}}
\def\xJ{{x_{\indJ}}}
\def\xE{x_{\indE}}
\def\po{\preceq}
\def\fI{f_{\indI}}
\def\fJ{f_{\indJ}}

\def\PI{P_{\indI}}
\def\PJ{P_{\indJ}}
\def\PE{P_{\indE}}
\def\phiI{\phi_{\indI}}
\def\phiJ{\phi_{\indJ}}

\def\topV{\mathcal{T}_{\ind{V}}}
\def\borelV{\borel_{\ind{V}}}
\def\topws{\mathcal{T}_{w^{\ast}}}
\def\borelws{\borel_{w^{\ast}}}
\def\wstar{weak$^{\ast}$ }

\def\Parts{\indset}
\def\part{\mathcal{H}}
\def\simp{\triangle}
\def\simpI{\triangle_{\indI}}
\def\simpJ{\triangle_{\indJ}}

\def\mean{\mathbb{E}}
\def\XI{X_{\indI}}
\def\XJ{X_{\indJ}}
\def\XE{X_{\indE}}
\def\meanI{\mathbb{E}_{P_{\indI}}}
\def\meanJ{\mathbb{E}_{P_{\indJ}}}
\def\meanE{\mathbb{E}_{P_{\indE}}}
\def\bU{{\mathcal{U}}}
\def\QU{\mathcal{Q}(\bU)}

\def\MV{\mathbf{M}(V)}
\def\MQ{\mathbf{M}(\mathcal{Q})}

\def\CQ{\mathbf{C}(\mathcal{Q})}

\def\L0{\mathcal{L}_0}

\def\cont{\CQ}

\def\famE#1{{\bigl< #1 \bigr>}_{\indE}}

\def\pMeas{\mathbf{M}}
\def\binary#1{[#1]_2}
\def\igdensity{p_{\mbox{\tiny\rm IG}}}

\def\NIG{\mbox{\rm NIG}}
\def\Polya{P\'{o}lya }

\endlocaldefs

\begin{document}

\begin{frontmatter}

\title{Projective limit random probabilities on~Polish spaces}
\runtitle{Projective limit random probabilities}
\begin{aug}
  \author{\fnms{Peter} \snm{Orbanz}\corref{}\ead[label=e1]{p.orbanz@eng.cam.ac.uk}\thanksref{t1}}
  \thankstext{t1}{Research supported by EPSRC grant EP/F028628/1.}
  \address{Computational and Biological Learning Laboratory\\
     University of Cambridge\\
     \printead{e1}}
  \affiliation{University of Cambridge}
  \runauthor{P. Orbanz}
\end{aug}


\begin{abstract}
A pivotal problem in Bayesian nonparametrics is the construction
of prior distributions on the space $\mathbf{M}(V)$ of probability measures
on a given domain $V$.
In principle, such distributions on the infinite-dimensional space $\mathbf{M}(V)$ can
be constructed from their finite-dimensional marginals---the most
prominent example being the construction of the Dirichlet
process from finite-dimensional Dirichlet distributions.
This approach
is both intuitive and applicable to the construction of arbitrary
distributions on $\mathbf{M}(V)$, but also hamstrung by a number of technical difficulties.
We show how these difficulties can be resolved if the
domain $V$ is a Polish topological space, and give a representation
theorem directly applicable to the construction of any probability
distribution on $\mathbf{M}(V)$ whose first moment measure is well-defined.
The proof draws on a projective limit theorem of Bochner,
and on properties of set functions on Polish spaces
to establish countable additivity of the resulting random probabilities.
\end{abstract}



\received{\smonth{1} \syear{2011}}

\end{frontmatter}

\renewcommand{\labelitemi}{-}

\section{Introduction}
\label{sec:introduction}

A variety of ways exists to construct the Dirichlet process.
For this particular case of a random probability measure, the
spectrum of construction approaches ranges from the projective
limit construction from finite-dimensional Dirichlet distributions
proposed by \citet{Ferguson:1973} to the stick-breaking construction
of \citet{Sethuraman:1994}; see \eg the survey by
\citet{Walker:Damien:Laud:Smith:1999} for an overview.
Most of these constructions are bespoke representations
more or less specific to the Dirichlet. An
exception is the projective limit representation, which
can represent any probability distribution on the space
of probability measures.
However, several authors
\citep[e.g.][]{Ghosal:2010,Ghosh:Ramamoorthi:2002}
have noted technical problems arising
for this construction.
The key role of the Dirichlet process, and
the proven utility of its representation
by stick-breaking or by Poisson
processes, may account for the slightly
surprising fact that these problems have not
yet been addressed in the literature.

The purpose of this paper is to provide a projective limit
result directly applicable to the construction of any
probability distribution on $\MV$. We do so by first
modifying and then proving a construction idea put forth by
\citet{Ferguson:1973}. Intuitively speaking,
our main result (Theorem \ref{theorem:main}) allows us to
construct distributions on $\MV$ by substituting
the Dirichlet distributions used in the derivation of the
Dirichlet process by other families of
distributions, and by verifying that these families satisfy
the two necessary and sufficient conditions of the theorem.
Stick-breaking, urn schemes \citep{Blackwell:MacQueen:1973}
and other specialized representations of the Dirichlet process all rely on
the latter's particular discreteness and spatial decorrelation properties.
Our approach may facilitate the derivation of models for which no such
representations can be expected to exist, for example, of smooth random measures.
For Bayesian nonparametrics, the result
provides what currently seems to be the only available tool to construct an arbitrary
prior distribution on the set $\MV$.
It also makes Bayesian methods based on random measures
more readily comparable to other types of nonparametric priors
constructed in a similar fashion,
notably to Gaussian processes \cite{Bauer:1996,Zhao:2000}.

The technical difficulties arising for the construction
proposed in \citep{Ferguson:1973} can be
summarized as three separate problems,
which Appendix \ref{sec:problems} reviews in detail. In short:
\begin{enumerate}
\renewcommand{\labelenumi}{\roman{enumi}}
\item {\bf Product spaces.} The product space setting of the standard
  Kolmogorov extension theorem is not well-adapted to the problem of
  constructing random probability measures.
\item {\bf Measurability problems.} A straightforward formalization
  of the construction in terms of an extension or projective
  limit theorem results in a space whose dimensions are labeled
  by the Borel sets of $V$, and is hence of uncountable dimension.
  As a consequence, the constructed measure cannot resolve most events
  of interest. In particular, singletons, and hence the event that the
  random measure assumes a specific measure as its value, are not measurable
  \citep[][Sec.\ 2.3.2]{Ghosh:Ramamoorthi:2002}.
\item {\bf $\sigma$-additivity.} The constructed measure is
  supported on finitely additive probabilities (charges), rather than
  $\sigma$-additive probabilities (measures);
  see \citet[][Sec. 2.2]{Ghosal:2010}.
  Further conditions
  are necessary to obtain a measure on probability measures.
\end{enumerate}

To make the projective limit construction feasible, we
have to impose some topological requirements on the
domain $V$ of the random measure. Specifically, we require
that $V$ is a Polish space, \ie a topological space which
is complete, separable and metrizable
\citep{Kechris:1995}.
This setting is sufficiently
general to accommodate any applications in Bayesian
nonparametrics---Bayesian methods do not solicit the generality of arbitrary
measurable spaces, since no useful notion of conditional probability
can be defined without a modicum of topological structure.
Polish spaces are in many regards
the natural habitat of Bayesian statistics, whether parametric or
nonparametric, since they guarantee both the existence of regular
conditional probabilities and the validity of de Finetti's theorem
\citep[][Theorem 11.10]{Kallenberg:2001}. The restriction to Polish spaces
is hence unlikely to incur any loss of generality.
We address problem (i) by means of a generalization
of Kolmogorov's extension
theorem, due to \citet{Bochner:1955}; problem (ii) by
means of the fact that the Borel $\sigma$-algebra of
a Polish space $V$ is generated by a countable subsystem
of sets, which allows us to substitute the uncountable-dimensional
projective limit space by a countable-dimensional surrogate;
and problem (iii) using a result of \citet{Harris:1968} on
$\sigma$-additivity of set functions on Polish spaces.

The remainder of the article is structured as follows:
The main result is stated in Sec.~\ref{sec:introduction:result}, which
is meant to provide all information required to apply the theorem,
without going into the details of the proof.
Related work is summarized in
Sec.~\ref{sec:introduction:related:work}.
A brief overview of projective limit constructions is given
in Sec.~\ref{sec:projlim}, to the extent relevant to the proof.
Secs.~\ref{sec:simplices} and \ref{sec:additivity} contain the
actual proof of Theorem \ref{theorem:main}: The projective
limit construction of random set functions
is described in Sec.~\ref{sec:simplices}.
A necessary and sufficient condition for these random set functions
to be $\sigma$-additive is given in Sec.~\ref{sec:additivity}.
Appendix \ref{sec:problems} reviews problems (i)-(iii) above
in more detail.

\subsection{Main result}
\label{sec:introduction:result}

To state our main theorem, we must introduce some notation,
and specify the relevant notion of a marginal distribution in the present context.
Let $\MV$ be the set of Borel probability measures over a Polish topological
space $(V,\topV)$; recall that the space is Polish if $\topV$ is a metrizable
topology under which $V$ is complete and separable \citep{Aliprantis:Border:2006,Kechris:1995}.
Throughout, the underlying model of randomness is an abstract probability space
$(\abstspace,\abstfield,\abstmeasure)$.
A random variable $X\!\!:\abstspace\rightarrow\MV$, with the image measure $P\!:=X\abstmeasure$ as its distribution,
is called a {\em random probability measure} on $V$.
Our main result, Theorem \ref{theorem:main}, is a general representation
result for the distribution $P$ of such a random measure.
To define measures on the space $\MV$, we endow it with the \wstar
topology $\topws$ (which in the context of probability is often called the topology of weak convergence)
and with the corresponding Borel $\sigma$-algebra
$\borelws\!:=\sigma(\topws)$. Since $V$ is Polish, the topological
space $(\MV,\topws)$ is Polish as well \citep[][Theorem 17.23]{Kechris:1995}.

Let $I=(A_1,\dots,A_n)$ be a \emph{measurable partition} of $V$,
\ie a partition of $V$ into a finite number of measurable, disjoint
sets. Denote the set of all such partitions $\part(\borelV)$.
Any probability measure $x\in\MV$ can be evaluated on a
partition $I$ to produce a vector $x_{\indI}:=(x(A_1),\dots,x(A_n))$,
and we write ${\phiI:x\mapsto x_{\indI}}$ for the evaluation functional so defined.
Clearly, $x_{\indI}$ represents a probability measure on the finite
$\sigma$-algebra $\sigma(I)$ generated by the partition.
Let $\simpI$ be the set of all
measures $x_{\indI}=\phiI(x)$
obtained in this manner, where $x$ runs through all measures
in $\MV$. This set, $\simpI=\phiI\MV$, is precisely the unit simplex in
the $n$-dimensional Euclidean space $\mathbb{R}^{\indI}$, 
\begin{equation}
 \label{eq:def:xspaceI}
  \simpI = \Bigl\lbrace \xI \in \mathbb{R}^{\indI}\Bigl\vert \xI(A_i)\geq 0
  \text{ and } \sum_{A_i\in I} \xI(A_i) = 1 \Bigr.\Bigr\rbrace \;.
\end{equation}

Let $J=(B_1,\dots,B_m)$ and $I=(A_1,\dots,A_n)$ be partitions such
that $I$ is a coarsening of $J$, that is, for each $A_i\in I$,
there is a set $\mathcal{J}_i\subset\lbrace 1,\dots,m\rbrace$
of indices such that $A_i=\cup_{j\in\mathcal{J}_i}B_j$. The sets
$\mathcal{J}_i$ form a partition of the index set
$\lbrace 1,\dots,m\rbrace$. If $I$ is a coarsening of $J$,
we write $I\po J$.

Let $x,x'\in\MV$. If $I\po J$, then $\phiJ x=\phiJ x'$ implies
$\phiI x=\phiI x'$. In other words, $\phiI x$ is completely determined
by $\phiJ x$, and invariant under any changes to $x$ which do not
affect $\phiJ x$. Therefore, the implicit definition
$\fJI(\phiJ(x)):=\phiI(x)$ determines a well-defined mapping
$\fJI:\simpJ\rightarrow\simpI$.
With notation for $J$ and $I$ as above, $\fJI$ can equivalently
be defined as
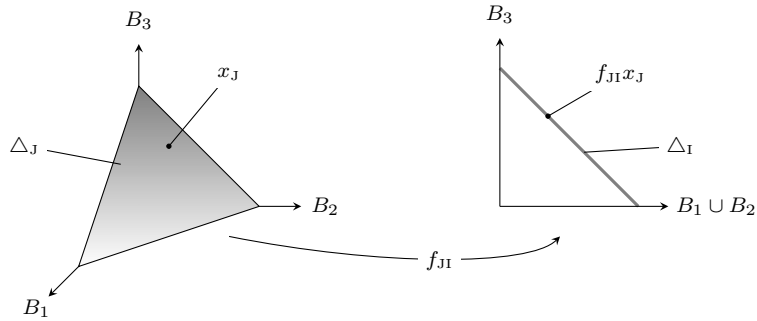
\begin{figure}
\begin{tikzpicture}[>=stealth,scale=0.8]
  \begin{scope}
    \draw[->] (0,0) -- (2.7,0); \draw (3.1,0) node {$B_2$};
    \draw[->] (0,0) -- (0,2.7); \draw (0,3.1) node {$B_3$};
    \draw[->] (0,0) -- (-1.5,-1.5); \draw (-1.7,-1.7) node {$B_1$};
    \shadedraw (2,0) --(0,2) --(-1,-1) --(2,0);
    \draw (-1.9,1) node [fill=white] {$\simpJ$} --(-0.3,0.7);
    \draw (1.5,2.2) node [fill=white] {$\xJ$} --(0.5,1);
    \filldraw (0.5,1) circle (1pt);
  \end{scope}
  \begin{scope}[xshift=6cm]
    {\color{black!50!white}
    \shadedraw[very thick] (2.3,0) --(0,2.3);}
    \draw[->] (0,0) -- (2.8,0); \draw (3.6,0) node {$B_1\cup B_2$};
    \draw[->] (0,0) -- (0,2.8); \draw (0,3.2) node {$B_3$};
    \draw (3,1) node [fill=white] {$\simpI$} --(1.4,0.9);
    \draw (2,2.2) node [fill=white] {$\fJI\xJ$} --(0.8,1.5);
    \filldraw (0.8,1.5) circle (1pt);
  \end{scope}
  \draw[->] (1.5,-0.5) .. controls (4,-1) and (6.5,-1) .. (7,-0.5);
  \draw (5,-0.9) node [fill=white] {$\fJI$};
\end{tikzpicture}
  \caption{Left: The simplex $\simpJ\subset\mathbb{R}^{\indJ}$ for a partition $J=(B_1,B_2,B_3)$.
    Right: A new simplex $\simpI=\fJI\simpJ$ is obtained by merging the sets $B_1$ and
    $B_2$, producing the partition $I=(B_1\cup B_2,B_3)$. The mapping $\fJI$ is given by
    $\fJI\xJ = (\xJ(B_1)+\xJ(B_2),\xJ(B_3))$.
    Its image $\simpI$
    is a subset of the product space $\mathbb{R}^{\indI}$, which shares only a single axis,
    $B_3$, with the space $\mathbb{R}^{\indJ}$.
  }
  \label{fig:simplex}
\end{figure}
\begin{equation}
        (\fJI\xJ)(A_i)=\sum_{j\in\mathcal{J}_i}\xJ(B_j) \;.
\end{equation}
Figure \ref{fig:simplex} illustrates
the mapping $\fJI$ and the simplices $\simpJ$ and $\simpI$.
The image $\fJI\xJ\in\simpI$ constitutes a probability distribution
on the events in $I$.
The following intuition is often helpful:
The space $\MV$ is convex, with the
Dirac measures on $V$ as its extreme points, and we can roughly think
of $\MV$ as the infinite-dimensional analogue of the simplices
$\simpI$. Similarly, we can regard the evaluations maps
$\phiI:\MV\rightarrow\simpI$
as analogues of the maps $\fJI:\simpJ\rightarrow\simpI$.
Even though both $\MV$ and the spaces $\simpI$ are Polish,
however, we have to keep in mind that the \wstar topology
on $\MV$ is, in many regards, quite different from the
topology which $\simpI$ inherits from Euclidean space.
For further properties of the space $\MV$, we
refer to the excellent exposition given by
\citet[][Chapter 15]{Aliprantis:Border:2006}.

Suppose that $P$ is a probability measure on $\MV$. Denote by $\phiI P$ the
image measure of $P$ under $\phiI$, \ie the measure on $\simpI$ defined
by $(\phiI P)(A_{\indI}):=P(\phiI^{-1}A_{\indI})$ for all $A_{\indI}\in\borel(\simpI)$.
We refer to $\phiI P$ as the \emph{marginal} of $P$ on $\simpI$. Similarly,
if $\PJ$ is a measure on $\simpJ$, then for any $I\po J$, the image measure
$\fJI\PJ$ is called the marginal of $\PJ$ on $\simpI$. The following theorem,
our main result, states that a measure $P$ on $\MV$ can be constructed from
a suitable family of marginals $\PI$ on the simplices $\simpI$.
The notation $\mean_Q[\,.\,]$ refers to expectation with respect to the law $Q$.

\newpage
\begin{theorem}
  \label{theorem:main}
  Let $V$ be a Polish space with Borel sets $\borelV$. Let
  $\MV$ be the set of probability measures on $(V,\borelV)$, and
  $\borelws$ the Borel $\sigma$-algebra generated by the \wstar topology
  on $\MV$. Let
  $\left<\PI\right>_{\part(\borelV)}:=\lbrace\PI|I\in\part(\borelV)\rbrace$
  be a family
  of probability measures on the finite-dimensional simplices
  $\simpI$. The following statements are equivalent:\\
  (1) The family $\left<\PI\right>_{\part(\borelV)}$ is projective,
  \begin{equation}
    \label{eq:intro:projective}
    \PI = \fJI\PJ \qquad\qquad\text{ whenever }I\po J
  \end{equation}
  and satisfies
  \begin{equation}
    \label{eq:intro:mean}
    \mean_{\PI}[\XI] = \phiI G_0 \qquad\qquad\text{ for all }I\in\part(\borelV)\;.
  \end{equation}
  (2)   There exists a unique probability measure $P$ on
  $(\MV,\borelws)$ satisfying
  \begin{equation}
    \label{eq:intro:marginals}
    \PI = \phiI P \qquad\qquad\text{ for all }I\in\part(\borelV)
  \end{equation}
  and
  \begin{equation}
    \mean_P[X]=G_0 \qquad\qquad\text{ for some }G_0\in\MV \;.
  \end{equation}
  If either statement holds, $P$ is a Radon measure.
\end{theorem}

\begin{remark}
  Theorem \ref{theorem:main} is applicable to the construction of
  any random probability measure $X$ on $V$ whose first moment
  $\mean_{P}[X]$ exists.
  In particular, the random measure $X$ need not be discrete.
  See Sec.~\ref{sec:examples} for examples.
\end{remark}

The two conditions of Theorem
\ref{theorem:main} serve two separate purposes:
Condition \eqref{eq:intro:projective} guarantees
that the family $\left<\PI\right>_{\part(\borelV)}$
defines a
unique probability measure $P_{\part(\borelV)}$. The support of this measure
is not actually $\MV$,
but a larger set---specifically, the set
$\CQ$ of
{\em finitely} additive probability measures (charges)
defined on a certain subsystem $Q\subset\borelV$, which we
will make precise in Sec.~\ref{sec:simplices}.
The set $\CQ$ contains the set $\MQ$ of $\sigma$-additive probability
measures on $Q$ as a measurable subset, and $\MQ$ is in turn
isomorphic to $\MV$, by Carath\'{e}odory's extension theorem \citep[][Theorem 2.5]{Kallenberg:2001}.
To obtain the distribution of a random measure, we need to
ensure that $P_{\part(\borelV)}$ concentrates on the subset $\MQ\cong\MV$, or
in other words, that draws from $P_{\part(\borelV)}$ are $\sigma$-additive
almost surely.
Condition \eqref{eq:intro:mean} is sufficient---and in
fact necessary---for $P_{\part(\borelV)}$ to concentrate on $\MV$, and
therefore for a random variable $X_{\part(\borelV)}$ with distribution $P_{\part(\borelV)}$
to constitute a random measure.
If \eqref{eq:intro:mean} is satisfied, the measure constructed
on $\CQ$ can be restricted to a measure on $\MV$,
resulting in the measure $P$ described by Theorem \ref{theorem:main}.
Sec.~\ref{sec:simplices} provides more details.

The technical restriction that $V$ be Polish is a mild one for all practical
purposes, a fact best illustrated
by some concrete examples of Polish spaces:
The real line is Polish, and so are $\mathbb{R}^n$ and $\mathbb{C}^n$;
any finite space; all separable Banach spaces (since Banach spaces are
complete metric spaces), in particular $\mathcal{L}_2$ and
any other separable Hilbert space; the space $\MV$ of probability measures
over a Polish domain $V$, in the \wstar topology \citep[][Chapter 15]{Aliprantis:Border:2006}; the
spaces $\mathcal{C}([0,1],\mathbb{R})$ and $\mathcal{C}(\mathbb{R}_+,\mathbb{R})$ of continuous functions, in the
topology of compact convergence \citep[][\S 38]{Bauer:1996}; and the Skorohod space
$\mathcal{D}(\mathbb{R}_+,\mathbb{R})$ of c\`{a}dl\`{a}g functions \citep[Chapter VI]{Pollard:1984}.
Any countable product of Polish spaces is Polish, in particular $\mathbb{R}^{\mathbb{N}}$,
$\mathbb{C}^{\mathbb{N}}$, and the Hilbert cube $[0,1]^{\mathbb{N}}$.
A subset of a given Polish space is Polish in the relative topology if and only
if it is a $G_{\delta}$ set \citep[][Theorem 3.11]{Kechris:1995}.
A borderline example are the spaces $\mathcal{C}(T,E)$ of continuous functions
with Polish range $E$. This space is Polish if $T=\mathbb{R}_{\geq 0}$
or if $T$ is compact and Polish, but not \eg for $T=\mathbb{R}$ \citep[][\S 454]{Fremlin}.
In Bayesian nonparametrics, this distinction may be relevant in the context of
the ``dependent Dirichlet process'' model of \citet{MacEachern:2000}, which involves
Dirichlet processes on spaces of continuous functions. For more background on
Polish spaces, see
\citep{Aliprantis:Border:2006,Fremlin,Kechris:1995}.

\subsection{Examples}
\label{sec:examples}

Theorem \ref{theorem:main} yields straightforward constructions for several models
studied in the literature, and we consider three specific examples
to illustrate the result.
First, by choosing the finite-dimensional marginals $\PI$ in
Theorem \ref{theorem:main} as a suitable family of
Dirichlet distributions, we obtain a construction
of the Dirichlet process in the spirit
of Ferguson \cite{Ferguson:1973}.
\begin{corollary}[Dirichlet Process]
  \label{corollary:DP}
  Let $V$ be a Polish space, $G_0$ a probability measure
  on $\borelV$, and let $\alpha\in\mathbb{R}_{>0}$.
  For each $I\in\part(\borelV)$, define $\PI$ as the Dirichlet distribution
  on $\simpI\subset\mathbb{R}^{\indI}$, with concentration $\alpha$
  and expectation $\phiI G_0\in\simpI$.  Then there is a uniquely determined
  probability measure $P$ on $\MV$ with expectation
  $G_0$ and the distributions $\PI$ as its marginals, that is,
  $\phiI P = \PI$
  for all $I\in\part(\borelV)$.
\end{corollary}

A similar construction yields the \emph{normalized inverse Gaussian process} of \citet{Lijoi:Mena:Pruenster:2005}.
The inverse Gaussian distribution on $\mathbb{R}_{\geq 0}$ is given by the density
$\igdensity(z|\alpha,\gamma)=\frac{\alpha}{\sqrt{2\pi}}x^{-3/2}\exp(-\frac{1}{2}(\frac{\alpha^2}{x}+\gamma^2 x)+\gamma\alpha)$ with respect to Lebesgue measure.
\citet{Lijoi:Mena:Pruenster:2005} define a normalized inverse Gaussian distribution $\NIG(\alpha_1,\dots,\alpha_n)$
on the simplex $\simp_n\subset\mathbb{R}^n$ as the distribution of the vector $w=(\frac{z_1}{\sum_i z_i},\dots,\frac{z_n}{\sum_i z_i})$,
where $z_i$ is distributed according to $\igdensity(z_i|\alpha_i,\gamma=1)$. The density of $w$ can be derived explicitly
\cite[][Equation (4)]{Lijoi:Mena:Pruenster:2005}.
Applicability of Theorem \ref{theorem:main} is a direct consequence of the results of \citet{Lijoi:Mena:Pruenster:2005},
which imply conditions
\eqref{eq:intro:projective} \citep[][(C3)]{Lijoi:Mena:Pruenster:2005} and
\eqref{eq:intro:mean} \citep[][Proposition 2]{Lijoi:Mena:Pruenster:2005}.
\begin{corollary}[Normalized Inverse Gaussian Process]
  \label{corollary:NIG}
  Let $\alpha\in\mathbb{R}_+$ and $G_0\in\MV$. For any partition $I=(A_1,\dots,A_n)$ in $\part(\borelV)$,
  choose the measure $\PI$ as the normalized inverse Gaussian distribution
  $\NIG(\alpha G_0(A_1),\dots,\alpha G_0(A_n))$.
  There is a uniquely determined
  probability measure $P$ on $\MV$ with expectation
  $G_0$ and $\phiI P = \PI$ for all $I\in\part(\borelV)$.
\end{corollary}

Although both the Dirichlet process and the normalized inverse Gaussian process are discrete almost surely,
Theorem \ref{theorem:main} is applicable to the construction of continuous random measures.
The \emph{\Polya tree} random measures introduced by \citet{Ferguson:1974} provide a convenient example.
They can be obtained as projective limits as follows:
Choose $V=\mathbb{R}$ and let $G_0\in\pMeas(\mathbb{R})$ be a probability measure with cumulative distribution function $g_0$.
For each $n$, let $I_n$ be the partition of $\mathbb{R}$ into intervals
$[g_0^{-1}(\frac{k-1}{2^n}),g_0^{-1}(\frac{k}{2^n}))$, where $k=1,\dots,2^n$.
All sets in $I_n$ have identical probability $1/2^n$ under $G_0$.
Since each partition $I_n$ is obtained from
$I_{n-1}$ by splitting each set in $I_{n-1}$ at a single point, the sequence $(I_n)$
satisfies $I_1\po I_2\po\dots$. It can be represented as a binary tree whose $n$th level corresponds
to $I_n$, each node representing one constituent set. There are two
natural ways of indexing sets in the partitions: One is to write $A_{n,k}$ for the $k$th set
in $I_n$, \ie $n$ indexes tree levels and $k$ enumerates sets within each level.
The other is to index sets as $A_{m_1,\dots,m_n}$ by a binary sequence encoding the unique path
from the root node $\mathbb{R}$ and the set in question, where $m_i=1$ indicates
passing to a right child node. Let $\binary{m}$
denote the binary representation of an arbitrary positive integer $m$. Then
\begin{equation*}
  A_{n,k}=\Bigl[g_0^{-1}\Bigl(\frac{k-1}{2^n}\Bigr),g_0^{-1}\Bigl(\frac{k}{2^n}\Bigr)\Bigr)=A_{\binary{2^n+(k-1)}}
    \quad\text{and}\quad
    I_n=(A_{n,1},\dots,A_{n,2^n})\;.
\end{equation*}
It is useful to use both index conventions interchangeably.
With each node $A_{m_1\cdots m_n}$, we associate a pair
$(Y_{m_1\cdots m_n 0},Y_{m_1\cdots m_n 1})\sim\mbox{Beta}(\alpha_{m_1\cdots m_n 0},\alpha_{m_1\cdots m_n 1})$
of beta random variables:
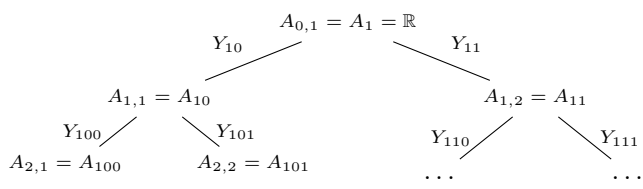
\begin{equation*}
\begin{tikzpicture}[baseline=(current bounding box.center),level/.style={sibling distance=50mm/#1},level distance=10mm]
  \begin{scope}[]
  \node [font=\scriptsize] (a00) {$A_{0,1}=A_1=\mathbb{R}$}
  child {node [font=\scriptsize] (a11) {$A_{1,1}=A_{10}$}
    child {node [font=\scriptsize] (a21) {
        \shortstack{
          $A_{2,1}=A_{100}$
          \\ \\
          $\dots$
        }
      }
    }
    child {node [font=\scriptsize] (a22)
      {
        \shortstack{
          $A_{2,2}=A_{101}$
          \\ \\
          $\dots$
        }
      }
    }
  }
  child {node [font=\scriptsize] (a12) {$A_{1,2}=A_{11}$}
    child {
      node (a23) {\shortstack{\mbox{\ }\\ $\dots$}}
    }
    child {
      node (a24) {\shortstack{\mbox{\ }\\ $\dots$}}
    }
  };
  \end{scope}
  \begin{scope}[font=\scriptsize]
  \path (a00) -- (a11) node[font=\scriptsize,scale=0.01,pos=0.5,label=above left:{$Y_{10}$}] {};
  \path (a00) -- (a12) node[font=\scriptsize,scale=0.01,pos=0.5,label=above right:{$Y_{11}$}] {};
  \path (a11) -- (a21) node[font=\scriptsize,scale=0.6,pos=0.5,label=left:{$Y_{100}$}] {};
  \path (a11) -- (a22) node[font=\scriptsize,scale=0.6,pos=0.5,label=right:{$Y_{101}$}] {};
  \path (a12) -- (a23) node[font=\scriptsize,scale=0.6,pos=0.5,label=left:{$Y_{110}$}] {};
  \path (a12) -- (a24) node[font=\scriptsize,scale=0.6,pos=0.5,label=right:{$Y_{111}$}] {};
  \end{scope}
\end{tikzpicture}
\end{equation*}

To apply Theorem \ref{theorem:main},
define probability measures $P_{\indI_{\ind{n}}}$ on the simplices $\simp_{\indI_{\ind{n}}}$ as follows:
Suppose a particle slides down the tree, moving along
each edge with the associated probability $Y_{m_1\cdots m_n}$.
The probability of reaching the set $A_{n,k}$ is a random variable $X_{n,k}$, defined recursively
in terms of the beta variables as $X_{m_1\cdots m_n m_{n+1}}\!\!:=X_{m_1\cdots m_n}Y_{m_1\cdots m_n m_{n+1}}$.
Choose $P_{\indI_{\ind{n}}}$ as the distribution of $X_{\indI_{\ind{n}}}=(X_{n,1},\dots,X_{n,2^n})$.
Applicability of Theorem \ref{theorem:main} follows from two results of
\citet{Ferguson:1974}:
(a) The partitions $I_n$ generate the Borel sets $\borel(\mathbb{R})$ and
(b) each random measure $X_{\indI_{\ind{n}}}\in\simp_{\indI_{\ind{n}}}$ has expectation
$\mean[X_{\indI_{\ind{n}}}]=(G_0(A_{n,1}),\dots,G_0(A_{n,2^n}))$.
Property (a) implies that the sequence
$P_{\indI_{\ind{n}}}$ induces a complete family $\left<\PI\right>$ of probability measures on
all simplices $\simp_{\indI}$, $I\in\part(\borel(\mathbb{R}))$.
By construction, $\left<\PI\right>$ satisfies \eqref{eq:intro:projective}. According to (b), \eqref{eq:intro:mean}
holds. Theorem \ref{theorem:main} and the well-known
continuity properties of \Polya trees \citep[][Theorem 3]{Lavine:1992} yield:
\begin{corollary}[\Polya tree]
  Let $\left<\PI\right>$ be a family of measures defined as above. There is a unique probability measure
  $P$ on $\pMeas(\mathbb{R})$ satisfying $\phiI P=\PI$. The distribution $P$ is a \Polya tree in the
  sense of \citet{Ferguson:1974}, with parameters $G_0$ and $(\alpha_{\binary{n}})_{n\in\mathbb{N}}$.
  The random probability measure
  $X$ on $\mathbb{R}$ with distribution $P$ has expected measure $\mean_P[X]=G_0$. If
  $\alpha_{n,k}=cn^2$ for some $c>0$, then $X$ is absolutely continuous with respect to Lebesgue measure
  on $\mathbb{R}$ almost surely.
\end{corollary}

\subsection{Related work}
\label{sec:introduction:related:work}

Theorem \ref{theorem:main} was effectively
conjectured by \citet{Ferguson:1973}. Although he only considered
the special case of the Dirichlet process, and despite the
technical difficulties already mentioned, he recognized both
the usefulness of indexing spaces by measurable partitions
(a key ingredient of the construction in Sec.~\ref{sec:simplices}),
and the connection between $\sigma$-additivity of
random draws from the Dirichlet process and $\sigma$-additivity
of its parameter measure
\citep[cf.][Proposition 2]{Ferguson:1973}.
Authors who have recognized problems to the effect that such
a construction is not feasible on an arbitrary measurable space
$V$ include \citet{Ghosh:Ramamoorthi:2002}
and \citet{Ghosal:2010}; both references also provide excellent surveys
of the different construction approaches available
for the Dirichlet process. \citet{Ghosal:2010} additionally points out,
in the context of problem (ii),
that a countable generator may be substituted for $\borelV$, provided
the underlying space is separable and metrizable.

To resolve the $\sigma$-additivity problem (iii), we appeal to a result
of \citet{Harris:1968}, which reduces the conditions for $\sigma$-additivity
of random set functions to their behavior on a countable number of sequences.
This result is well-known in the theory
of point processes and random measures
\cite{Kallenberg:1983,Crauel:2002}.
Although Sethuraman was aware of Harris' work and
referenced it in his well-known article
\citep{Sethuraman:1994}, it has to our knowledge never been
followed up on in the nonparametric Bayesian literature.

For the specific problem of defining the Dirichlet process, it is
possible to forego the projective limit construction altogether and
invoke approaches specifically tailored to the properties of the
Dirichlet
\citep{Ghosal:2010,Ghosh:Ramamoorthi:2002,Walker:Damien:Laud:Smith:1999}.
On the real line, both the Dirichlet process and the
closely related Poisson-Dirichlet distribution of Kingman
\citep{Kingman:1975} arise in a variety of contexts throughout
mathematics, each of which can be regarded as a possible means of
definition \citep[e.g.][]{Olshanski:2003,Talagrand:2003}.
On arbitrary Polish spaces, the Dirichlet process can be derived
implicitly as de Finetti mixing measure of an
urn scheme \citep{Blackwell:MacQueen:1973},
or as special case of a P\'{o}lya tree \cite{Ferguson:1974}.

Sethuraman's stick-breaking scheme \citep{Sethuraman:1994}
is remarkable not only for its simplicity. In contrast to all other
constructions listed above, it does not require $V$ to be Polish, but is
applicable on an arbitrary measurable space with measurable singletons.
The stick-breaking and projective
limit representations of the Dirichlet process
trade off two different types of generality:
Stick-breaking imposes less
restrictions on the choice of $V$, but is not applicable to
represent other types of distributions
on $\MV$. The projective limit approach requires more structure
on $V$, but can represent any probability measure on $\MV$.
The trade-off is
reminiscent of similar phenomena encountered throughout
stochastic process theory. For example, probability measures
on infinite-dimensional product spaces can be constructed
by means of Kolmogorov's extension theorem. If the measure
to be constructed is factorial over the product, the component
spaces of the product may be chosen as arbitrary measurable spaces
\citep[][Theorem 9.2]{Bauer:1996}.
To model stochastic
dependence across different subspaces, however, a minimum of topological
structure is indispensable, and Kolmogorov's theorem hence requires
the component spaces to be Polish \citep[][Theorem 6.16]{Kallenberg:2001}. The
Dirichlet process, as a purely atomic random measure whose
different atoms are stochastically dependent only through the
global normalization constraint, can be regarded as the closest
analogue of a factorial measure on the space $\MV$. In analogy
to a factorial measure, it can be constructed on very general
spaces, whereas the projective limit approach, which can represent
arbitrary correlation structure, requires stronger topological
properties.

\section{Background: Projective limits}
\label{sec:projlim}

A projective limit is constructed from a family of mathematical
structures, indexed by the elements of an index set $\indset$
\cite{Bourbaki:1966,Mallory:Sion:1971}. For our
purposes, the structures in question will be topological
measurable spaces $(\xspaceI,\borelI)$, with $I\in\indset$.
The projective limit defined by this family is again a
measurable space, denoted $(\xspaceE,\borelE)$. This projective
limit space is the smallest space containing
all spaces $(\xspaceI,\borelI)$ as its substructures,
in a sense to be made precise shortly.
To obtain a meaningful notion of a limit, the index set $\indset$
need not be totally ordered, but it must be possible to form infinite
sequences of suitably chosen elements. The set is therefore required to
be \emph{directed}: There is a partial order relation $\po$ on $\indset$
and, whenever $I,J\in\indset$, there exists $K\in\indset$
such that $I\po K$ and $J\po K$.
A simple example of a directed set is the set $\indset:=\mathcal{F}(L)$ of
all finite subsets of an infinite set $L$, where $\indset$ is partially
ordered by inclusion.

The component spaces $\xspaceI$ used to define the projective limit
need to ``fit in'' with each other in a suitable manner. This idea is
formalized by defining a family of mappings $\fJI$ between the
spaces which are regular with respect to the
structure posited on the point sets $\xspaceI$. For measurable spaces,
the adequate notion of regularity is measurability. Since we assume
each $\sigma$-algebra $\borelI$ to be generated by an underlying topology
$\topI$, we slightly strengthen this requirement to continuity.
\begin{definition}[Projective limit set]
  Let $\indset$ be a directed set and $(\xspaceI,\topI)$, with $I\in
  \indset$, a family of topological spaces.
  For any pair $I\po J\in\indset$, let
  $\fJI:\xspaceJ\rightarrow\xspaceI$ be a function such that
  \begin{enumerate}
  \item $\fJI$ is $\topJ$-$\topI$-continuous.
  \item $f_{\indI\indI}=\mbox{Id}_{\xspaceI}$.
  \item $\fKJ\circ\fJI=\fKI$ whenever $I\po J\po K$.
  \end{enumerate}
  The functions $\fJI$ are called \emph{generalized projections}.
  The family $\lbrace \xspaceI,\topI,\fJI | I\po J\in
  \indset\rbrace$, which we denote $\famE{\xspaceI,\topI,\fJI}$,
  is called a {\em projective system} of topological spaces.
  Define a set $\xspaceE$ as follows:
  For each collection $\lbrace \xI\in\xspaceI | I \in\indset\rbrace$
  of points satisfying
  \begin{equation}
    \label{eq:def:proj:lim:point}
    \xI=\fJI\xJ \qquad\qquad\text{ whenever } I\po J \;,
  \end{equation}
  identify the set $\lbrace \xI\in\xspaceI | I \in\indset\rbrace$
  with a point $\xE$, and let $\xspaceE$ be the collection of all
  such points. The set $\xspaceE$ is called the \emph{projective limit set}
  of $\famE{\xspaceI,\fJI}$.
\end{definition}
Denote the Borel
$\sigma$-algebras on the topological spaces $\xspaceI$ by $\borelI\!\!:=\sigma(\topI)$.
For each $I\in\indset$, the map defined as $\fI\!\!:\xE\mapsto\xI$ is a
well-defined function $\fI\!\!:\xspaceE\rightarrow\xspaceI$. These functions
are called \emph{canonical mappings}.
They define a topology $\topE$ and a $\sigma$-algebra on the projective
limit space $\xspaceE$, as the smallest topology (resp.\ $\sigma$-algebra)
which makes all canonical mappings $\fI$ continuous (resp.\ measurable).
In particular,
\begin{equation}
  \label{eq:projlim:sigma:algebra}
  \borelE:=\sigma( \fI | I\in\indset) = \sigma(\cup_{I\in\indset}
  \fI^{-1}\borelI) = \sigma(\topE)\;.
\end{equation}
In analogy to the set $\xspaceE$, the topological
space $(\xspaceE,\topE)$ is called the projective limit of $\famE{\xspaceI,\topI,\fJI}$,
and the measurable space $(\xspaceE,\borelE)$ the
projective limit of $\famE{\xspaceI, \borelI, \fJI}$.

A measure $\PE$ on the projective limit $(\xspaceE,\borelE)$ can
be constructed by defining a measure $\PI$ on each
space $(\xspaceI,\borelI)$. By simultanously applying the projective
limit to the projective system $\famE{\xspaceI,\borelI,\fJI}$ and
to the measures $\PI$, the family $\famE{\PI}$ is
assembled into the measure $\PE$. The only requirement is that
the measures $\PI$ satisfy a condition analogous to the one imposed
on points by \eqref{eq:def:proj:lim:point}.
More precisely, $\PI$ has to coincide with the image measure of $\PJ$
under $\fJI$,
\begin{equation}
  \label{eq:def:projective:family}
  \PI = \fJI\PJ = \PJ\circ\fJI^{-1} \qquad\qquad\text{ whenever
  }I\po J\;.
\end{equation}
A family of measures $\famE{\PI}$ satisfying
\eqref{eq:def:projective:family} is called a {\em projective family}.
The existence and uniqueness of $\PE$ on $(\xspaceE,\borelE)$ is
guaranteed by the following result
\citep[][IX.4.3, Theorem 2]{Bourbaki:2004}.

\begin{theorem}[Bochner]
  \label{theorem:bochner}
  Let $\famE{\xspaceI, \borelI, \fJI}$ be a
  projective system of measurable spaces with countable, directed index set $D$,
  and $\famE{\PI}$
  a projective family of probability measures on these spaces.
  Then there exists a uniquely defined measure $\PE$ on the
  projective limit space $(\xspaceE,\borelE)$ such that
  \begin{equation}
    \PI = \fI\PE = \PE\circ\fI^{-1} \qquad\qquad\text{ for all }I\in\indset\;.
  \end{equation}
\end{theorem}
We refer to the measures in the family $\famE{\PI}$ as the \emph{marginals} of
the stochastic process $\PE$. Since the marginals completely determine
$\PE$, some authors refer to $\famE{\PI}$
as the \emph{weak distribution}
of the process, or as a \emph{promeasure} \cite{Bourbaki:2004}.

Theorem \ref{theorem:bochner} was introduced by
\citet[][Theorem 5.1.1]{Bochner:1955}, for a possibly uncountable index set $\indset$.
The uncountable case requires an additional condition known
as \emph{sequential maximality}, which ensures the
projective limit space is non-empty. For our purposes, however, countability of
the index set is essential: Measurability problems
(problem (ii) in Sec.~\ref{sec:introduction})
arise whenever $D$ is uncountable, and are not resolved by sequential
maximality.

The most common example of a projective limit theorem in probability theory
is Kolmogorov's extension theorem \cite[][Theorem 6.16]{Kallenberg:2001}, which
can be regarded as the special case of Bochner's theorem obtained for product spaces:
Let $\indset$ be the set of all finite subsets of an infinite set $L$, partially
ordered by inclusion. Choose any Polish measurable space $(\xspace_0,\borel_0)$,
and set $\xspaceI:=\prod_{i\in I}\xspace_0$. The resulting projective limit space
is the infinite product $\xspaceE=\prod_{i\in L}\xspace_0$,
and $\borelE$ coincides with the Borel
$\sigma$-algebra generated by the product topology. For product spaces, the
sequential maximality condition mentioned above holds automatically, so $L$
may be either countable or uncountable. Once again, though, the measurability
problem (ii) arises unless $L$ is countable.
The product space form of the theorem is typically used in the
construction of Gaussian process distributions on random functions
\cite{Bauer:1996}. For
random measures, a more adequate projective system is constructed in
following section.

\section{Projective limits of probability simplices}
\label{sec:simplices}

This section constitutes the first part of the proof of Theorem
\ref{theorem:main}: The construction of a projective limit space
$\xspaceE$ from simplices $\simpI$, and the analysis of its properties.
The space $\xspaceE$ turns out to consist of set functions which
are not necessarily $\sigma$-additive, and the remaining part of
the proof in Sec.~\ref{sec:additivity} will be the derivation of
a criterion for $\sigma$-additivity.

The distinction between finitely additive and
$\sigma$-additive set functions will be crucial to the ensuing
discussion. We consider two types of set systems $\mathcal{Q}$
on the space $V$: {\em Algebras},
which contain both $\emptyset$ and $V$, and are closed under
complements and finite unions, and {\em $\sigma$-algebras}, which
are algebras and additionally closed under countable unions.
A non-negative set function $\mu$ on either an algebra or $\sigma$-algebra
$\mathcal{Q}$ is called a \emph{charge} if it satisfies
$\mu(\emptyset)=0$ and is finitely additive. If a charge
is normalized, \ie if $\mu(V)=1$, it is called a \emph{probability charge}.
A charge is a measure if and only if it is $\sigma$-additive. If $\mathcal{Q}$ is an algebra,
and not closed under countable unions,
the definition of $\sigma$-additivity only requires
$\mu$ to be additive along those countable
sequences of sets $A_n\in\mathcal{Q}$ whose union is in $\mathcal{Q}$.

\subsection{Definition of the projective system}
\label{sec:simplices:proj:system}

For the choice of components in a projective system, it can
be helpful to regard the elements $\xE$ of the projective limit
space $\xspaceE$ as mappings,
from a domain defined by the index set $\indset$ to a range defined
by the spaces $\xspaceI$.
The simplest example is once again the product space
$\xspaceE=\xspace_0^{\ind{L}}$
in Kolmogorov's theorem, for which
each $\xE\in\xspaceE$ can be interpreted as a function
$\xE:L\rightarrow\xspace_0$.
Probability measures on $(V,\borelV)$ are in particular set functions
$\borelV\rightarrow [0,1]$, so it is natural to construct $\indset$ from
the sets in $\borelV$. It is not necessary to include all measurable
sets: If $\mathcal{Q}$ is an algebra that generates $\borelV$,
any probability measure on $\mathcal{Q}$ has, by Carath\'{e}odory's theorem
\cite[][Theorem 2.5]{Kallenberg:2001}, a unique extension to a probability measure on $\borelV$.
In other words, the space $\MQ$ of probability measures on
$\mathcal{Q}$ is isomorphic to $\MV$, and $\mathcal{Q}$ can be substituted for
$\borelV$ in the projective limit construction.

Desiderata for the projective limit are:
(1) The projective limit space $\xspaceE$ should contain all measures
on $\mathcal{Q}$ (and hence on $\borelV$). (2) $\mathcal{Q}$ should
be countable, to address the measurability problem (ii) in Sec.~\ref{sec:introduction}.
(3) The marginal spaces $\xspaceI$ should consist of the finite-dimensional
analogues of measures on $\mathcal{Q}$, and hence of measures on
finite subsets of events in $\mathcal{Q}$. (4) The definition of
the system should facilitate a proof of $\sigma$-additivity.
In this section, we will recapitulate the projective limit specified
in Sec.~\ref{sec:introduction:result} and show it indeed satisfies
(1)-(3); that (4) is satisfied as well will be shown in
Sec.~\ref{sec:additivity}.
\\

{\noindent\bf Choice of $\mathcal{Q}$.}
We start with the prototypical choice of basis for any Polish topology:
Let $\denseset\subset V$ be a countable, dense subset of $V$. Fix
a metric ${d:V\times V\rightarrow \mathbb{R}_+}$ which
generates the topology $\topV$, and denote by $B(v,r)$
the open $d$-ball of radius $r$ around $v$. Denote the
system of open balls with rational radii and centers in $\denseset$ by
\begin{equation}
  \label{eq:def:bU}
  \bU := \lbrace B(v,r) | v\in\denseset, r\in\mathbb{Q}_+\rbrace \cup \lbrace\emptyset\rbrace \;.
\end{equation}
Since $V$ is separable and metrizable, $\mathcal{U}$ forms
a countable basis of the topology $\topV$
\citep[][Lemma 3.4]{Aliprantis:Border:2006}.
Let $\QU$ be the algebra generated by $\bU$. Then
$\bU\subset\QU\subset\borelV$. In particular,
$\QU$ is a countable generator of $\borelV$.
\\

{\noindent\bf Index set.}
As the index set $\indset$, we do not choose $\QU$ itself, but rather
the set of all finite partitions of $V$ consisting of
disjoint sets $A_i\in\QU$,
\begin{equation}
  \label{eq:def:Parts}
  \indset:=
  \part(\mathcal{Q})
  = \Bigl\lbrace (A_1,\dots,A_n) \,\Bigl\vert\, n\in\mathbb{N},\,
  A_i\in\QU,\, \dot{\cup}A_i = V \Bigr.\Bigr\rbrace \;.
\end{equation}
Each element $I\in D$ is a finite partition, and the set of
probability measures on the events in this partition is precisely
the simplex $\simpI$.
To define a partial order on $\Parts$, let
$I=(A_1,\dots,A_m)$ and $J=(B_1,\dots,B_n)$ be any two partitions in
$\Parts$, and denote their intersection (common refinement) by
$I\cap J := (A_i\cap B_j)_{i,j}$.
Since $\QU$ forms an algebra, $I\cap J$ is again an element of $\Parts$.
Now define a partial order relation $\po$ as
\begin{equation}
  I \po J \quad :\Leftrightarrow \quad I\cap J=J \;,
\end{equation}
that is, $I\po J$ iff $J$ is a refinement of $I$. The set $(\Parts,\po)$
is a valid index set for a projective limit system, because it is directed:
$K:=I\cap J$ always satisfies $I\po K$ and
$J\po K$.
\\

{\noindent\bf Projection functions.}
What remains to be done is to specify the functions $\fJI$.
Consider a partition $J=(A_1,\dots,A_n)$, and any
$\xJ\in\simpJ$. Each entry $\xJ(A_j)$ assigns a number (a
probability) to the event $A_j$, and we define $\fJI$ accordingly
to preserve this property.
To this end, let $J=(B_1,\dots,B_n)$ be a partition in $\Parts$, and let
$I=(A_1,\dots,A_m)$ be a coarsening of $J$ (that is, $I\po J$).
For each $A_i$, let $\mathcal{J}_i\subset\lbrace 1,\dots,n\rbrace$
be the subset of indices for which
$A_i = \cup_{j\in\mathcal{J}_i} B_j$.
Then define $\fJI$ as
\begin{equation}
  \label{eq:def:fJI}
  (\fJI\xJ)(A_i) := \sum_{j\in\mathcal{J}_i} \xJ(B_j) \;.
\end{equation}
We choose $\xspaceI:=\simpI$ as defined in \eqref{eq:def:xspaceI},
and endow $\simpI$ with the relative
topology $\topI:=\top(\mathbb{R}^{\indI})\cap\simpI$ and the
corresponding Borel sets
${\borelI\!:=\borel(\topI)=\borel(\mathbb{R}^{\indI})\cap\simpI}$.
The relative topology makes additions on $\simpI$, and
hence the mappings $\fJI$, continuous. Each $f_{\indI\indI}$ is
the identity on $\simpI$, and $f_{\indK\indI}=\fKJ\circ\fJI$.
For any pair $I\po J\in\Parts$, $\simpI=\fJI\simpJ$ and
conversely, $\simpJ=\fJI^{-1}\simpI$.
Therefore, $\famE{\simpI,\borelI, \fJI}$
is a projective system.

\subsection{Structure of the projective limit space}

Let $(\xspaceE,\borelE)$ be the projective limit of
$\famE{\simpI,\borelI,\fJI}$.
We observe immediately that $\xspaceE$ contains $\MQ$:
If $x$ is a probability measure on $\QU$,
let $\xI:=\fI x$ for each partition $I\in\Parts$.
The collection $\lbrace \xI | I\in\Parts\rbrace$ satisfies
\eqref{eq:def:proj:lim:point}, and hence constitutes a point in $\xspaceE$.
The following result provides more details about the
constructed measurable space
$(\xspaceE,\borelE)$, which turns out to
be the space $\CQ$ of all probabiliy charges defined on $\QU$.
By $\borelws$, we again denote the Borel $\sigma$-algebra
on $\MV$ generated by the \wstar topology.
\begin{proposition}
  \label{theorem:structure:xspaceE}
  Let $V$ be a Polish space, and $(\xspaceE,\borelE)$ the projective limit
  of the projective system $\famE{\simpI,\borelI,\fJI}$ defined in
  Sec.~\ref{sec:simplices:proj:system}.
  Denote by $\psi:\MV\rightarrow\MQ$ the restriction mapping
  which takes each measure $x$ on $\borelV$ to its restriction
  $\xE=x\vert_{\mathcal{Q}}$ on $\mathcal{Q}\subset\borelV$.
  Then the following hold:
  \begin{enumerate}
    \renewcommand{\labelenumi}{(\roman{enumi})}
  \item $\xspaceE=\CQ$, the space of probability charges on $\QU$.
  \item $\MQ$ is a measurable subset of $\CQ$.
  \item $\psi$ is a Borel isomorphism of $(\MV,\borelws)$ and
    $(\MQ,\borelE\cap\MQ)$.
  \end{enumerate}
\end{proposition}
Part $(ii)$ implies that
a projective limit measure $\PE$ constructed on $\CQ$
by means of Theorem \ref{theorem:bochner} can
be restricted to a measure on $\MQ$ without further complications,
in particular without appealing to outer measures.
According to $(iii)$, there is a measure $P$ on $\MV$ which
can be regarded as equivalent to $\PE$, namely the image
measure $P:=\psi^{-1}\PE$ under the inverse of the restriction
map $\psi$. This is of course the measure $P$ described in
Theorem \ref{theorem:main}, though some details still remain
to be established later on. Since $\psi$ is a Borel isomorphism,
$P$ constitutes a measure with respect to the ``natural''
topology on $\MV$.

\begin{proof}

  {\noindent\em Part (i).} Let $\xE\in\xspaceE$.
  The trivial partition $I_0:=(V)$ is in $\Parts$, which implies
  $\xE(V)= f_{I_0}\xE=1$ and $\xE(\emptyset)=0$. To show finite
  additivity, let $A_1,A_2\in\QU$ be disjoint sets and choose a
  partition $J\in\Parts$ such that $A_1,A_2\in J$. Let $I\po J$ be the
  coarsening of $J$ obtained by joining the two sets.
  As the elements of each space $\simpI$ are finitely additive,
  \begin{equation}
    \xE(A_1)+\xE(A_2)=(f_{\indJ}\xE)(A_1)+(f_{\indJ}\xE)(A_2)
    \stackrel{\eqref{eq:def:fJI}}{=}
    (f_{\indI}\xE)(A_1\cup A_2)=\xE(A_1\cup A_2) \;.\nonumber
  \end{equation}
  Hence, $\xE$ is a charge.
  Conversely, assume that $\xE$ is a probability charge on $\QU$.
  The evaluation $\fI\xE$ of $\xE$ on a partition $I\in\Parts$
  defines a probability measure on the finite $\sigma$-algebra
  $\sigma(I)$, and thus $\fI\xE\in\simpI$.
  Since additionally $\fJI(\fJ\xE)=\fI\xE$, the set
  $\famE{\fI\xE}$ forms a collection
  of points $\fI\xE\in\simpI$ satisfying
  \eqref{eq:def:proj:lim:point}, and hence $\xE\in\xspaceE$.

  \noindent\emph{Part (ii).}
  Regard the restriction map $\psi$ as a mapping
  into $\CQ$, with image $\MQ$. By Caratheodory's extension theorem,
  $\psi$ is injective \citep[][Theorem 2.5]{Kallenberg:2001}.
  If an injective mapping between Polish spaces is measurable,
  its inverse is measurable as well \citep[][Theorem A1.3]{Kallenberg:2001}.
  Thus, if we can show $\psi$ to be measurable, $\MQ=\psi(\MV)$
  is a measurable set.

  First observe that $\psi$ relates the
  evaluation functionals ${\fI:\CQ\rightarrow\simpI}$ on
  probability charges to the evaluation functionals $\phiI:\MV\rightarrow\simpI$ on
  probability measures via the equations
  \begin{equation}
    \label{eq:relation:phiI:fI:psi}
    \phiI = \fI\circ\psi \qquad\text{ for all } I\in D\;.
  \end{equation}
  We will show that the mappings $\phiI$ generate the $\sigma$-algebra $\borelws$ on
  $\MV$. Since the canonical mappings $\fI$ generate $\borelE$ on $\CQ$ by definition,
  \eqref{eq:relation:phiI:fI:psi} then implies $\borelws$-$\borelE$-measurability of $\psi$:

  Let $\phi_A\!:\MV\rightarrow[0,1]$ be the evaluation functional $x\mapsto x(A)$.
  Since $\MV$ is separable, the Borel sets of the \wstar topology coincide with those
  generated by the maps $\phi_A$ \citep[][Theorem 2.3]{Gaudard:Hadwin:1989}, thus
  $\borelws=\sigma( \phi_A | A\in\borelV)$. Each mapping $\phi_A$ can be identified
  with $\phiI$ for $I=(A,A^c)$, because $\phi_{(A,A^c)}(x)=(x(A),1-x(A))$.
  Hence equivalently, $\borelws=\sigma( \phi_{(A,A^c)} | A\in\borelV)$, and with
  \eqref{eq:relation:phiI:fI:psi},
  \begin{equation}
    \borelws=\psi^{-1}\sigma( f_{(A,A^c)} | A\in\borelV)\;.
  \end{equation}
  Clearly, the maps $f_{(A,A^c)}$ for $A\in\mathcal{Q}$ are sufficient
  to express all information expressible by the larger family of maps
  $\fI$, $I\in D$, and thus generate the projective limit $\sigma$-algebra,
  \begin{equation}
    \sigma(f_{(A,A^c)} \,\vert\, A\in\mathcal{Q})=\borelE \;.
  \end{equation}
  In summary, $\psi$ is $\borelws$-$\borelE$-measurable, and we deduce $\MQ=\psi(\MV)\in\borelE$.

  \noindent\emph{Part (iii).} As shown above, $\psi$ is injective and
  measurable, and regarded as a mapping onto its image $\MQ$, it is trivially
  surjective. What remains to be shown is measurability of the inverse.
  By part $(ii)$, the image $\psi(\MV)=\MQ$ is a Borel subset of $\CQ$.
  As a countable projective limit of Polish spaces, $(\CQ,\topE)$ is Polish
  \citep[][Chapter IX]{Bourbaki:1966}.
  Since $\MV$ is Polish, $(\MV,\borelws)$ is a standard Borel space, \ie a
  Borel space generated by a Polish topology.
  The space $(\MQ,\borelE\cap\MQ)$ is standard Borel as well, since
  $\MQ$ is a Borel subset of a Polish space \citep[][Theorem A1.2]{Kallenberg:2001}.
  As noted above, measurable bijections between standard Borel spaces are automatically
  bimeasurable \citep[][Theorem A1.3]{Kallenberg:2001}, which shows
  $\psi$ to be a Borel isomorphism.
\end{proof}

\section[$sigma$-additivity of random charges]{$\sigma$-additivity of random charges}
\label{sec:additivity}

The previous section provides the means to construct
the distribution $\PE$ of a random charge $\XE:\abstspace\rightarrow\cont$
as a projective limit measure. To obtain random measures
rather than random charges in this manner, we need to additionally
ensure that $\PE$ concentrates on the measurable subspace $\MV$,
or in other words, that $\XE$ is $\sigma$-additive $\abstmeasure$-almost surely.

Consider a projective limit random charge $\XE$, distributed
according to a projective limit measure $\PE$ on $\cont$.
The following proposition gives a necessary and sufficient condition
for almost sure $\sigma$-additivity of $\XE$,
formulated in terms of its expectation $\mean_{\PE}[\XE]$.
It also shows that the expected values of $\PE$ and the projective
family $\famE{\PI}$ are themselves projective, in the sense
that $\fI\mean_{\PE}[\XE]=\mean_{\PI}[\XI]$, and accordingly
$\fJI\mean_{\PJ}[\XJ]=\mean_{\PI}[\XI]$ for any pair $I\po J$.
The latter makes the criterion directly applicable to construction
problems: If we initiate the construction by choosing an expected
measure $G_0\in\MV$ for the prospective measure $\PE$, and then choose
the projective family such that $\mean_{\PI}[\XI]=\fI G_0$, random
draws from $\PE$ will take values in $\MV$ almost surely.

\begin{proposition}
  \label{theorem:mean:sigma:additivity}
  Let $(\xspaceE,\borelE)$ be the projective limit of
  finite-dimensional probability
  simplices defined in Proposition \ref{theorem:structure:xspaceE},
  and let $\famE{\PI}$ be a
  projective family of probability measures on the spaces
  $(\simpI,\borelI)$. Denote by $\PE$ the projective limit measure,
  and by $G_0:=\meanE[\XE]$ its expectation.
  Then:
  \begin{enumerate}
    \renewcommand{\labelenumi}{(\roman{enumi})}
    \item The expectation $G_0$ is an element of $\xspaceE$ and
      \begin{equation}
        \label{eq:marginals:projective}
        \fI G_0 = \meanI[\XI]\qquad\qquad\text{ for any } I\in\Parts \;.
      \end{equation}
    \item
      $\XE$ is $\sigma$-additive $\abstmeasure$-almost surely if and only if
      $G_0$ is $\sigma$-additive.
  \end{enumerate}
\end{proposition}

The proof requires a criterion for $\sigma$-additivity of probability charges
expressible in terms of a countable number of conditions. Assuming that
$G_0$ is $\sigma$-additive, we will deduce from
the projective limit construction that, if a fixed sequence of sets is given,
the random content $\XE$ is countably additive along this sequence with probability
one. This only implies almost sure $\sigma$-additivity of $\XE$ on $\QU$ if
the condition for $\sigma$-additivity can be reduced to a countable subset
of sequences in $\QU$ (cf.\ Appendix \ref{sec:problems:additivity}).
Such a reduction was derived by Harris \citep[][Lemma 6.1]{Harris:1968}.
For our particular choice of $\QU$, his result can be stated as follows:
\begin{lemma}[Harris]
  \label{lemma:sigma:additivity}
  Let $V$ be any Polish space and $\QU$ the countable algebra generated by the
  open balls \eqref{eq:def:bU}. Then the set of all sequences of elements of
  $\QU$ contains a countable subset of sequences
  $(A_n^m)_{n}$, where $A^m_n\searrow\emptyset$ for all $m\in\mathbb{N}$,
  such that any probability charge $\mu$ on $\QU$
  is $\sigma$-additive if and only if it satisfies
  \begin{equation}
    \label{eq:s-additivity:condition}
    \lim_{n\rightarrow\infty} \mu(A_n^m) = 0 \qquad\qquad\text{ for all } m\in\mathbb{N}\;.
  \end{equation}
\end{lemma}

\begin{proof}[Proof of Proposition \ref{theorem:mean:sigma:additivity}]

  {\noindent\em Part (i).}
  The expectation $\meanE[\XE]$ is finitely additive: For any finite number of disjoint
  sets $A_i\in\borelE$,
  \begin{equation}
    \sum_{i=1}^n\meanE[\XE](A_i)
    =
    \int_{\CQ}\sum_{i=1}^n\xE(A_i)\PE(d\xE)
    =
    \meanE[\XE](\cup_i A_i)\;.
  \end{equation}
  Since clearly also $\meanE[\XE](\emptyset)=0$ and
  $\meanE[\XE](V)=1$, the expectation is an element of $\xspaceE$.
  To verify \eqref{eq:marginals:projective}, note
  the mappings $\fJI:\simpJ\rightarrow\simpI$ are affine, and
  hence
  \begin{equation}
    \begin{split}
    \fJI\meanJ[\XJ]
    \stackrel{f \text{affine}}{=}
    \meanJ[\fJI\XJ]
    =&
    \int_{\simpJ=\fJI^{-1}\simpI}
    \fJI\xJ\PJ(d\xJ)\\
    =&
    \int_{\simpI}\xI(\fJI\PJ)(d\xI)
    =
    \int_{\simpI}\xI\PI(d\xI) \;.
    \end{split}
  \end{equation}
  Therefore, the expectations of a projective family $\famE{\PI}$
  satisfy $\fJI\meanJ[\XJ]=\meanI[\XI]$.
  By the same device,
  $\fI G_0 = \fI\meanE[\XE]=\meanI[\XI]$ holds for the projective limit measure $\PE$.

  {\noindent\em Part (ii).}
  First assume that $G_0$ is $\sigma$-additive.
  Let $(A_n^m)_n$ be any of the set sequences given by Lemma
  \ref{lemma:sigma:additivity}.
  As $n\rightarrow\infty$, the random sequence $(\XE(A_n^m))$ converges to $0$ almost surely:
  $\sigma$-Additivity of $G_0$ implies
  \begin{equation}
    \label{eq:closure:under:expectation:MQ}
    \lim_{n\rightarrow\infty} \meanE[\XE](A_n^m)
    =
    \lim_{n\rightarrow\infty} G_0(A_n^m)
    =
    G_0(\emptyset)
    =
    \meanE[\XE](\emptyset) \;,
  \end{equation}
  hence $\XE(A_n^m)\xrightarrow{L_1}0$.
  The sequence $(A_n^m)$ is decreasing and the random variable $\XE$ is charge-valued, which implies
  $\XE(A_{n+1}^m)\leq\XE(A_n^m)$ a.s. In particular, the sequence $(\XE(A_n^m))$ forms
  a supermartingale when endowed with its canonical filtration.
  For supermartingales, convergence in the mean implies almost sure
  convergence \citep[][Theorem 19.3]{Bauer:1996}, and thus indeed
  $\XE(A_n^m)\xrightarrow{\mbox{\tiny a.s.}}0$.

  Consequently, there is a $\abstmeasure$-null subset $N_m$ of the abstract
  probability space $\abstspace$ such that
  \begin{equation}
    (\XE(\omega))(A_n^m)
    \xrightarrow{n\rightarrow\infty}
    (\XE(\omega))(\emptyset)
    \qquad\qquad\text{ for } \omega\not\in N_m\;.
  \end{equation}
  The union $N:=\cup_{m\in\mathbb{N}}N_m$ of these null sets, taken over all sequences $(A_n^m)$
  required by Lemma \ref{lemma:sigma:additivity}, is again a
  $\abstmeasure$-null set. The charge $\XE(\omega)$
  satisfies \eqref{eq:s-additivity:condition}
  for all $m$ whenever $\omega\not\in
  N$. Therefore, $\XE$ is $\sigma$-additive $\abstmeasure$-a.s.
  by Lemma \ref{lemma:sigma:additivity}, and hence
  almost surely a probability measure.

  Conversely, let $\XE$ assume values in $\MV\cong \MQ$ almost surely.
  Since $A_n^m\searrow\emptyset$, the sequence of measurable functions $\omega\mapsto(\XE(\omega))(A_n^m)$
  converges to $0$ almost everywhere. By hypothesis, $\CQ\!\smallsetminus\!\MQ$ is a null set, hence
  \begin{equation}
    \label{eq:s-additivity:of:mean:convergence}
    \lim_{n\rightarrow\infty}
    \meanE[\XE](A_n^m)
    =
    \lim_{n\rightarrow\infty}
    \int_{\XE^{-1}\MQ}(\XE(\omega))(A_n^m)\abstmeasure(d\omega)
    =
    0 \;,
  \end{equation}
  where the second identity holds
  by dominated convergence \citep[][Theorem 1.21]{Kallenberg:2001}. Since
  $\meanE[\XE]$ is a probability charge according to part $(i)$ and satisfies
  \eqref{eq:s-additivity:of:mean:convergence}, it satisfies the conditions of Lemma \ref{lemma:sigma:additivity},
  and we conclude $\meanE[\XE]\in\MQ$.
\end{proof}

Theorem \ref{theorem:main} is now finally obtained by deducing the properties of $P$
from those of $\PE$ as established by Proposition
\ref{theorem:mean:sigma:additivity}.

\begin{proof}[Proof of Theorem \ref{theorem:main}]
  First suppose that \eqref{eq:intro:projective} and \eqref{eq:intro:mean} hold.
  By Theorem \ref{theorem:bochner}, a unique projective limit measure $\PE$
  exists on $\xspaceE=\CQ$, with $\fI\PE=\PI$. Proposition
  \ref{theorem:mean:sigma:additivity}$(ii)$ shows $\PE$ is concentrated
  on the measurable subset
  $\MQ$. By Proposition \ref{theorem:structure:xspaceE}$(iii)$,
  it uniquely defines an equivalent measure $P:=\psi^{-1}\PE$ on $\MV$,
  which satisfies \eqref{eq:intro:marginals}. As a probability measure on a
  Polish space, $P$ is a Radon measure \citep[IX.3.3, Proposition 3]{Bourbaki:2004}.

  Conversely, assume that $P$ is given. Then \eqref{eq:intro:projective} follows from
  \eqref{eq:intro:marginals}.
  The expectation $G_0=\mean_{P}[X]$ is in $\MV$ by Proposition \ref{theorem:mean:sigma:additivity}$(ii)$.
  Any measure on $\MV$ can
  be represented as a measure on $\xspaceE=\CQ$, hence by
  Proposition \ref{theorem:mean:sigma:additivity}$(i)$, the expectation
  $G_0$ and the marginals $\PI=\fI P$ satisfy \eqref{eq:marginals:projective}.
  Thus, \eqref{eq:intro:mean} holds, and the proof is complete.
\end{proof}

\appendix
\section{Review of technical problems}
\label{sec:problems}

This appendix provides a more detailed description of
problems (i)--(iii) listed in Sec.~\ref{sec:introduction}.
The discussion addresses readers of passing familiarity
with measure-theoretic probability; to the probabilist,
it will only state the obvious.

The approach proposed in \citep{Ferguson:1973} is, in summary,
the following: A probability measure on $(V,\borelV)$ is
a set function $\borelV\rightarrow[0,1]$.
The set $\MV$ of probability measures can be regarded as a
subset of the space $[0,1]^{\borelV}$ of all such functions.
More precisely, the space chosen in \citep{Ferguson:1973}
is $[0,1]^{\part(\borelV)}$, where $\part(\borelV)$ again denotes
the set of all measurable, finite partitions of $V$.
This space contains one axis for each partition, and hence
is a larger space than $[0,1]^{\borelV}$, but
redundantly encodes the same information.
The Kolmogorov extension theorem
\citep[][Theorem 6.16]{Kallenberg:2001} is then applied
to a family of Dirichlet distributions defined on the
finite-dimensional subspaces of the product space
$[0,1]^{\part(\borelV)}$.

\subsection{Product spaces}
The Kolmogorov extension theorem used in the construction
is not well-adapted to the problem of constructing measures on
measures, because the setting assumed by the theorem is that
of a product space:
A finite-dimensional marginal of a measure $P$ on $\MV$
is a measure $\PI$ on the set of measures over a finite
$\sigma$-algebra $\mathcal{C}$ of events. Any such $\sigma$-algebra
can be generated by a partition $I$ of events in $\borelV$.
The set consisting of the marginals on $I$ of all measures $x\in\MV$
is necessarily isomorphic to the unit simplex in $|I|$-dimensional
Euclidean space. Hence, the marginals of a measure $P$ defined \emph{on}
$\MV$ always
live on simplices of the form $\simpI$ as described in
Sec.~\ref{sec:introduction:result}.
In other words, when we set up a projective limit construction
for measures on $\MV$, the choice of possible finite-dimensional
marginal spaces is limited---either the simplices are used directly,
as in Sec.~\ref{sec:introduction:result}, or they are embedded into
some other finite-dimensional space.
If the projective limit result to be applied is the Kolmogorov extension
theorem, the simplices must be embedded into Euclidean product
spaces, as proposed in \citep{Ferguson:1973}.
The problem here is that it is difficult to properly formalize
marginalization to subspaces, as required by the theorem.
For constructions on $[0,1]^{\borelV}$, the problem can be
illustrated by the example in Fig.~\ref{fig:simplex}:
For $J=(B_1,B_2,B_3)$, the simplex $\simpJ$ is a subspace
of $\mathbb{R}^{\indJ}\cong\mathbb{R}^3$.
Marginalization corresponds to merging two events,
such as $B_1$ and $B_2$ in the example. The resulting
simplex $\simpI$ for $I=(B_1\cup B_2,B_3)$ is a subspace
of $\mathbb{R}^{\indI}$. However, $\mathbb{R}^{\indI}$
is not a subspace of $\mathbb{R}^{\indJ}$, nor is
$\simpI$ a subspace of $\simpJ$.
Hence, in the product space setting of the Kolmogorov theorem,
the natural way to formalize a reduction in dimension
for measures on a finite number of events
does not correspond to a projection onto a subspace.

\subsection{Measurability problems}

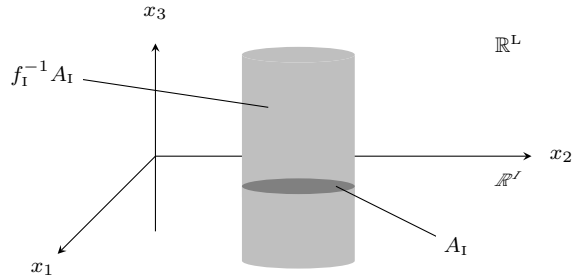
\begin{figure}
\begin{tikzpicture}[>=stealth]
  \draw[->] (0,0) -- (5,0); \draw (5.4,0) node {$x_2$};
  \draw[->] (0,-1) -- (0,1.5); \draw (0,1.9) node {$x_3$};
  \draw[->] (0,0) -- (-1.3,-1.3); \draw (-1.5,-1.5) node {$x_1$};
  \node [scale=3, cylinder, stroke, rotate=90, cylinder uses custom fill, cylinder end fill=black!25,
    cylinder body fill=black!25, aspect=0.35, minimum height=0.7cm, minimum width=0.5cm] at (1.9,0.3) {} ;
  \node [scale=3, cylinder, stroke, rotate=90, cylinder uses custom fill, cylinder end fill=black!50,
    cylinder body fill=black!25, aspect=0.35, minimum height=0.4cm, minimum width=0.5cm] at (1.9,-1) {} ;
  \draw (4,-1.2) --(2.4,-0.4); \draw (4,-1.2) node [fill=white] {$A_{\indI}$};
  \draw (-1.5,1.1) --(1.5,0.65); \draw (-1.5,1.1) node [fill=white] {$\fI^{-1}A_{\indI}$};
  \draw  (4.7,1.5) node [fill=white] {$\mathbb{R}^{\ind{L}}$};
  \draw (4.7,-0.3) node [fill=white,scale=0.85,xslant=0.7] {$\mathbb{R}^{\indI}$};
\end{tikzpicture}
  \caption{Three-dimensional analogue of a cylinder set in the product space setting.
    An event $A_{\indE}\subset\mathbb{R}^{\indE}$ is independent of the
    random variable $X_3$ if it is the preimage $A_{\indE}=\fI^{-1}A_{\indI}$ of some event
    $A_{\indI}\subset\mathbb{R}^{\indI}$, that is, if the set $A_{\indE}$ is of ``axis parallel'' shape
    in direction of $X_3$. The event $A_{\indI}$ in the figure occurs if $(X_1,X_2)\in A_{\indI}$, or
    equivalently, if $(X_1,X_2,X_3)\in\fJI^{-1}A_{\indI}$.}
  \label{fig:cylinder}
\end{figure}
A general property of projective limit constructions of
stochastic processes is that the index set---intuitively,
the set of axes labels of a product, or of dimensions in
a more general setting---must be countable to obtain
a useful probability measure. This is due to the fact
that all projective limit theorems implicitly generate a
$\sigma$-algebra on the infinite-dimensional space---the
$\sigma$-algebra $\borelE$
specified by \eqref{eq:projlim:sigma:algebra}---based
on the $\sigma$-algebras on the marginal spaces used
in the construction.
The constructed measure lives on this $\sigma$-algebra.

If the dimension is uncountable, the
resolution of the $\sigma$-algebra is too coarse to
resolve most events of interest. In particular, it does
not contain singletons.
The problem is most readily illustrated
in the product space setting: Suppose the Kolmogorov
theorem is used to define a measure $P$ on an infinite-dimensional
product space $\xspaceE\!:=\mathbb{R}^{\ind{L}}$, where $L$ is some infinite
set. The measure $P$ is constructed from given measures $\PI$
defined on the finite-dimensional sub-products $\mathbb{R}^{\indI}$,
where $I\in\indset$ are finite subsets of $L$. The $\sigma$-algebra
on $\mathbb{R}^{\ind{L}}$ on which $\PE$ is defined is generated
as follows: Denote by $\fI$ the product space projector
$\mathbb{R}^{\ind{L}}\rightarrow\mathbb{R}^{\indI}$. For any measurable
set $A^{\indI}\in\mathbb{R}^{\indI}$, the preimage $\fI^{-1}A_{\indI}$
is a subset of $\mathbb{R}^{\ind{L}}$, which is of ``axis-parallel''
shape in direction of all axis not contained in $I$. The finite-dimensional
analogue of this situation is illustrated in Fig.~\ref{fig:cylinder},
where $A^{\indI}$ is assumed to be an elliptically shaped set in the
plane $\mathbb{R}^{\indI}$, and the overall space $\mathbb{R}^{\ind{L}}$
is depicted as three-dimensional. Preimages $\fI^{-1}A_{\indI}$ of measurable
sets are, for obvious reasons, called \emph{cylinder sets} in the
probability literature.
The $\sigma$-algebra defined by the Kolmogorov theorem is the
smallest $\sigma$-algebra containing all cylinder sets $\fI^{-1}A_{\indI}$,
for all measurable sets $A_{\indI}\in\mathbb{R}^{\indI}$ and all
finite sub-products $\mathbb{R}^{\indI}$. Since $\sigma$-algebras
are defined by closure under countable operations, the sets in
this $\sigma$-algebra can be thought of as cylinder sets that
are of axis-parallel shape along all but a countable number of
dimensions. If the overall space is of countable dimension,
any set of interest can be expressed in this form. If the
dimension is uncountable, however, these events only specify
the joint behavior of a countable subset of random
variables---in Fig.~\ref{fig:cylinder}, $\mathbb{R}^{\indI}$ would represent
a subspace of countable dimension of the uncountable-dimensional
space $\mathbb{R}^{\ind{L}}$.

For example, consider the set
$\mathbb{R}^{\ind{L}}:=\mathbb{R}^{\mathbb{R}}$, regarded as the
set of all functions $\xE:\mathbb{R}\rightarrow\mathbb{R}$,
which arises in the construction of Gaussian processes.
Although the constructed measure $\PE$ is a distribution on random functions
$\xE$, this measure cannot assign a probability to events
of the form $\lbrace \XE=\xE\rbrace$, \ie to the event that the outcome
of a random draw is a particular function $\xE$. The only measurable
events are of the form $\lbrace \XE(s_1)=t_1,\XE(s_2)=t_2,\dots\rbrace$
and specify the value of the function at a countable subset of
points $s_1,s_2,\ldots \in\mathbb{R}$.

\subsection[$sigma$-additivity]{$\sigma$-additivity}
\label{sec:problems:additivity}

The marginal distributions used in the construction specify
the joint behavior of the constructed measure $\PE$ on any
finite subset of measurable sets. $\sigma$-additivity
requires additivity along an infinite sequence, and cannot
be deduced directly from additivity of the marginals.
Suppose that some sequence $A_1,A_2,\dots$ of measurable
sets in $V$ is given, and that $\xE$ is a random set function
drawn from $\PE$.
Countable additivity of $\xE$ along the sequence can
be shown to hold almost surely (with respect to $\PE$)
by means of a simple convergence argument
\cite[][Proposition 2]{Ferguson:1973}. However,
as a $\sigma$-algebra, $\borelV$ is either finite or uncountable.
Hence, if $V$ is infinite, $\borelV$ contains an uncountable
number of such sequences. Even though $\xE$ is additive
along any given sequence with probability one, the
null sets of exceptions aggregate into a non-null set
over all sequences, and $\xE$ is not $\sigma$-additive
with probability one. Substituting a countable generator
$\mathcal{Q}$ for $\borelV$ does not resolve the problem,
since the number of sequences in $\mathcal{Q}$ remains uncountable.

\section*{Acknowledgments}
I would like to thank the associate editor and two referees for valuable suggestions,
in particular for pointing out the example in Corollary \ref{corollary:NIG}.
I am grateful to Daniel M.\ Roy for helpful
comments and corrections.

\bibliographystyle{natbib}

\end{document}